\documentclass[leqno,11pt]{article}
\usepackage{latexsym}
\usepackage{amssymb}
\usepackage{amsfonts}
\usepackage{amsmath}
\usepackage{epsfig}
\usepackage{color}
\usepackage{pifont}

\usepackage{todonotes}

\allowdisplaybreaks
\makeatletter
\long\def\unmarkedfootnote#1{{\long\def\@makefntext##1{##1}\footnotetext{#1}}}
\makeatother

\newtheorem{definition}{Definition}[section]

\newtheorem{lemma}[definition]{Lemma}

\newtheorem{theorem}[definition]{Theorem}

\newtheorem{proposition}[definition]{Proposition}

\newtheorem{corollary}[definition]{Corollary}

\newtheorem{remark}[definition]{Remark}

\def\o{\Omega}

\def\t{\tau}

\def\m2{|\Omega | /2}
\def\M2{\frac{|\Omega |}{2}}
\def\u+{u_+^*}

\def\-p{\overline{p}}

\def\w0{{W_0^{1,p}(\Omega)}}

\def\MR{\mathcal R}
\def\R{\mathbb R}
\def\N{\mathbb N}

\def\ep{\varepsilon}

\def\rn{{{\R}^n}}
\def\T{T_{t}}

\newcommand{\hh}{{\cal H}^{n-1}}
\newcommand{\medint}{-\kern  -,395cm\int}
\newcommand{\medintinrigo}{-\kern  -,315cm\int}
\newcommand{\medelle}{-\kern  -,235cm L}
\newcommand{\medellenrigo}{-\kern  -,180cm L}
\newcommand{\qed}{\thinspace\null\nobreak\hfill
\hbox{\vbox{\kern-.2pt\hrule height.2pt
depth.2pt\kern-.2pt\kern-.2pt \hbox to1.8mm {\kern-.2pt\vrule
width.4pt \kern-.2pt\raise1.8mm\vbox to.2pt{} \lower0pt\vtop
to.2pt{}\hfil\kern-.2pt \vrule
width.4pt\kern-.2pt}\kern-.2pt\kern-.2pt \hrule height.2pt
depth.2pt \kern-.2pt}}\par\medbreak}

\title{
Second-order regularity for  parabolic $p$-Laplace problems
} 
\numberwithin{equation}{section}

\author{
  Andrea Cianchi\\
 {\it Dipartimento di Matematica e Informatica \lq\lq U. Dini", Universit\`a di Firenze}\\ {\it Viale Morgagni 67/A, 50134 Firenze, Italy} \\{\it  e-mail: cianchi@unifi.it}
\bigskip
\\
  Vladimir G. Maz'ya \\
  {\it   Department of Mathematics, Link\"oping University, SE-581
83 Link\"oping, Sweden}
  \\ and \\
{\it  RUDN University}\\ {\it
6 Miklukho-Maklay St, Moscow, 117198, Russia}
\\ {\it e-mail: vladimir.mazya@liu.se}
}

\date{}

\begin{document}
\maketitle
\begin{abstract}

Optimal second-order regularity  in the space variables is established  for solutions to Cauchy-Dirichlet problems for nonlinear parabolic equations and systems of $p$-Laplacian type, with 
square-integrable right-hand sides and  initial data in a Sobolev space. As a consequence, generalized solutions are shown  to be  strong solutions.
Minimal regularity  on the boundary of the domain is required, though the results are new even for smooth domains. In particular, they  hold in arbitrary bounded convex domains. 
%

\end{abstract}

\unmarkedfootnote {

%
%
}

\section{Introduction}\label{intro}



We deal with Cauchy-Dirichlet problems for parabolic equations and systems  of the form
\begin{equation}\label{eqdirichlet}
\begin{cases}
u_t - {\rm div} (|\nabla u|^{p-2} \nabla u ) = f  & {\rm in}\,\, \o_T\\
 u =0  &
{\rm on}\,\, \partial \o \times (0, T) \,
\\ u(\cdot ,0)= \psi (\cdot) &
\hbox{in} \,\, \Omega\,.
\end{cases}
\end{equation}
 Here, {\color{black} $p>1$,} $\Omega$ is an open set  in $\rn$, $n \geq 1$, with  finite Lebesgue measure $|\o|$,  and $T>0$. Moreover,
$$\Omega _T = \o \times (0, T),$$
 the functions $f: \o_T \to \mathbb R^N$ and $\psi : \o \to \mathbb R^N$, $N\geq 1$ are given, and   $u: \o \to \mathbb R^N$ is the unknown. According to usage, $\nabla u$ stands for the gradient of $u$ with respect to the space variables $x \in \o$, and $u_t$ for its  derivative in time $t \in (0,T)$.
\par We are concerned with global second-order  regularity properties, with respect to the variables $x$, of the solution $u$ to problem \eqref{eqdirichlet}. Our main results can be summarized as follows. Assume that
\begin{equation}\label{hpf}
f \in L^2(\o_T)\,,
\end{equation}
\begin{equation}\label{hppsi}
\psi \in W^{1,p}_0(\o)\,,
\end{equation}
and $\o$ satisfies suitable regularity conditions (if $n \geq 2$). Then
\begin{equation}\label{loose}
|\nabla u|^{p-2} \nabla u \in L^2((0,T); W^{1,2}(\o))\,,
\end{equation}
and the norm of $|\nabla u|^{p-2} \nabla u$ in $L^2((0,T); W^{1,2}(\o))$ is bounded by the norms of the data $f$ and $\psi$. Consequently, $u$ is actually a strong solution to problem \eqref{eqdirichlet}. 
%
This 
provides a    natural nonlinear counterpart of the classical $L^2((0,T); W^{2,2}(\o))$ regularity of solutions to Cauchy-Dirichlet problems for the heat equation \cite{Eid, Fr, LaSoUr, Li}.
\par The results of the present paper are  new even in the case of smooth domains $\o$.  They will however be established under minimal regularity assumptions on $\partial \o$. In particular, they hold in any convex bounded domain $\o$.
\par
Let us notice that, if $p\geq 2$, then assumptions \eqref{hpf} and \eqref{hppsi} ensure that $f$ and $\psi$ belong to  proper function spaces for a classical weak solution to problem \eqref{eqdirichlet} to be well defined. On the other hand, this is not guaranteed if $1<p<2$. Some specification is thus in order.
\par
When $N=1$, namely when dealing with a single equation, a generalized notion of solution can still be introduced to cover the whole range of exponents $p \in (1, \infty)$, and our regularity theory holds   for every such $p$. One kind of solution that fits the situation at hand can be defined as the limit of solutions to approximating problems involving smooth data \cite{BG, DallA, Pr}. Such a solution will be called approximable throughout this paper. Its existence, uniqueness and basic regularity under \eqref{hpf} and \eqref{hppsi} is established  in Theorem \ref{exist}.  Its second-order differentiability properties -- the central issue of our contribution --  are addressed  in Theorems \ref{secondcap} and \ref{secondconvex}.
\par When $N>1$, i.e. when systems are in question, we restrict our attention to the case $p\geq 2$. This is the subject of Theorems \ref{secondcapsys} and \ref{secondconvexsys}. A reason for this limitation on $p$ is the lack of a suitable  existence and uniqueness theory of approximable solutions for parabolic systems.
One obstacle for this gap is the failure of standard truncation methods for vector-valued functions. In fact, also the elliptic theory of approximable solutions  for systems is incomplete. This affects 
 our approach, which makes critical use of parallel results for elliptic equations and systems  recently obtained in \cite{cmsecond} and \cite{cmsecondsys}. 
\par To conclude this section, let us point out that, in spite of  huge developments of the regularity theory of nonlinear singular and degenarate parabolic problems, presented e.g. in the reference monographs \cite{Di, DiGiVe, Li, lions} and in recent papers including \cite{AKN, BaHa, BoDuMa, BoDuMi, DSSV, FS, KuM1, KuM2, Schw},
 information available in the literature about second-order regularity of solutions to nonlinear parabolic problems is still limited. A result in the spirit of \eqref{loose} can be found in  \cite{Beirao}, which yet  requires additional regularity of the datum $f$   and is restricted to values of $p$ smaller than and close to $2$. 
More classical contributions instead concern differentiability properties of the nonlinear expression $|\nabla u|^{\frac{p-2}2} \nabla u$, which differes from 
 that appearing in \eqref{loose}. Furthermore,  they tipically apply to local solutions, and request higher regularity of $f$.

\section{Main results}\label{main}

Let $N \geq 1$ and let $1<p<\infty$. Assume that $f\in L^{p'}((0,T); W^{-1,p'}(\o))$ and $\psi \in L^2(\o)$, where $p'= \tfrac p{p-1}$, the H\"older's conjugate of $p$. Then there exists a unique weak solution $u$ to problem \eqref{eqdirichlet}, namely a function $u \in C([0,T]; L^2(\o))\cap L^p((0,T); W^{1,p}_0(\o))$  such that $u_t \in L^{p'}((0,T); W^{-1,p'}(\o))$, $u(\cdot , 0) = \psi (\cdot)$ and
\begin{equation}\label{lions}
\int _0^\tau\langle u_t, \phi\rangle \, dt + \int _0^\tau \int_\o |\nabla u|^{p-2} \nabla u \cdot \nabla \phi \,dx\, dt  = \int _0^\tau\langle f, \phi\rangle \, dt  
\end{equation}
for every $\tau \in (0, T]$ and
 every  function $\phi \in L^p((0,T); W^{1,p}_0(\o))$. Here, $\langle \cdot , \cdot \rangle$ stands for the duality pairing between $W^{-1,p'}(\o)$ and $W^{1,p}_0(\o)$. This follows from classical monotonicity arguments -- see e.g. \cite[Chapter 2, Sections 1.1 and 1.5]{lions}. 

\par As mentioned above, if the assumptions that $f \in  L^{p'}((0,T); W^{-1,p'}(\o))$ and $\psi \in L^2(\o)$ are dropped, then the notion of weak solution to problem \eqref{eqdirichlet} does not apply  anymore. Still, in the spirit of  \cite{BG, DallA, Pr}, generalized solutions in the approximable sense can be defined whenever $f$ and $\psi$ are merely integrable, namely if  $f\in L^1(\o_T)$ and $\psi \in L^1(\o)$.
\\
A function $u \in C([0,T]; L^1(\o)) \cap L^1((0,T); W^{1,1}_0(\o))$   will be called an approximable solution to problem \eqref{eqdirichlet}  if there exist  sequences $\{f_k\}\subset C_0^\infty(\o_T)$  and $\{\psi_k\}\subset C^\infty_0(\o)$  such that
$$_k \to f \quad \hbox{in $L^1(\o_T)$,} \quad \quad \psi_k \to \psi \quad \hbox{in $L^1(\o)$,}$$
%
and the sequence $\{u_k\}$ of 
 weak solutions to the problems
\begin{equation}\label{eqdirichletk}
\begin{cases}
(u_k)_t - {\rm div} (|\nabla u_k|^{p-2} \nabla u_k ) = f_k  & {\rm in}\,\,\o_T\\
 u_k=0  &
{\rm on}\,\,
 \partial \o \times (0, T) \,
\\ u_k(\cdot ,0)= \psi_k (\cdot) &
{\rm in}\,\, \Omega\,
\end{cases}
\end{equation}
satisfies 
\begin{equation}\label{june0}
u_k \to u \quad \hbox{and} \quad \nabla u_k \to \nabla u \quad \hbox{a.e. in $\o_T$,}
\end{equation}
and 
\begin{equation}\label{june0bis}
u_k \to u \quad \hbox{in $ C([0,T]; L^1(\o))$.}
\end{equation}

Our first result ensures that, if $N=1$, problem \eqref{eqdirichlet} actually admits a unique approximable solution for every $p \in (1, \infty)$ which, under assumptions \eqref{hpf} and \eqref{hppsi}, enjoys additional regularity properties. Its uniqueness easily implies  that such a solution agrees with the weak solution if $p\geq 2$. This follows, for instance, by an argument as in the proof of Theorem \ref{secondcapsys}.

\begin{theorem}\label{exist} {\rm {\bf [Existence, uniqueness and basic regularity of approximable solutions]}}
Let $N=1$, $1<p<\infty$  and  $T>0$. Assume that $\o$ is an open set with finite measure in $\rn$,  \  $n \geq 1$. Let $f \in
L^2(\Omega_T )$ and $\psi \in W^{1,p}_0(\o)$. Then
there exists a unique  approximable solution $u$ to problem \eqref{eqdirichlet}.  Moreover,
\begin{equation}\label{solamain1}
 u\in  L^\infty((0,T); W^{1,p}_0(\o)) \quad \hbox{and}    \quad u_t\in L^2(\o_T)\,,
\end{equation}
and there exists a constant $C=C(n, p, |\o|)$ such that
\begin{align}\label{solamain} 
\| u\|_{L^\infty((0,T); W^{1,p}_0(\o))}^{\frac p2}+ \|u_t\|_{L^2(\o_T)}\leq C \big(\|f\|_{L^2(\o_T)} + \|\nabla \psi\|_{L^p(\o)}^{\frac p2}\big)\,.
\end{align}
%
%
%
%
\end{theorem}

Let us now turn to the primary objective of this paper, namely the global square integrability of weak derivatives  in $x$ of the  nonlinear expression of the gradient $|\nabla u|^{p-2} \nabla u$. 
\\ We begin by introducing a few notions to be used in prescribing the regularity of $\Omega$ when $n \geq 2$.
A minimal regularity condition on $\o$ for our result to hold involves   a local isocapacitary inequality for 
 the integral of the weak curvatures of $\partial \o$.  Specifically,
$\Omega$ is assumed to be a bounded Lipschitz domain. Moreover,  the functions of    $(n-1)$ variables that locally describe the boundary of $\o$ are  required to be twice weakly differentiable and integrable on $\partial \o$ with respect to the $(n-1)$-dimensional Hausdorff measure $\hh$. This will be denoted by  $\partial \o \in W^{2,1}$. In particular, the weak second fundamental form 
$\mathcal B$  on $\partial \o$ belongs to $L^1(\partial \o)$
These assumptions are not yet sufficient,
 as  demonstrated, for instance, by examples avaliable in the literature in the stationary  case -- see Remark \ref{sharp} below.  
A local smallness condition on the $L^1$-norm of $\mathcal B$ is also needed. Specifically, denote 
 by
$|\mathcal B|$ the norm of $\mathcal B$, and set 
\begin{equation}\label{defK}
\mathcal K_\o(r) =
\displaystyle 
\sup _{
\begin{tiny}
 \begin{array}{c}{E\subset \partial \o \cap B_r(x)} \\
x\in \partial \o
 \end{array}
  \end{tiny}
} \frac{\int _E |\mathcal B|d\hh}{{\rm cap}_{B_1(x)} (E)}\qquad \hbox{for $r\in (0,1)$}\,.
\end{equation}
Here, $B_r(x)$ stands for the ball centered at $x$, with radius $r$,  and the notation ${\rm cap}_{B_1(x)}(E)$ is adopted for the capacity of the set $E$ relative to the ball $B_1(x)$.
Then we require that
\begin{equation}\label{capcond}
\lim _{r\to 0^+} \mathcal K _\o(r) < c_1
\end{equation}
for a suitable constant $c_1=c_1(n, N, p,d_\o, L_\o)$,
where, 
  $d_\o$ and $L_\o$  denote the diameter
  and the   Lipschitz constant of $\o$, respectively. Here, and in similar occurrences in what follows, the dependence of a constant on $d_\o$  and  $L_\o$ is understood just via an upper bound for them. 
\\ Recall that 
the capacity ${\rm cap}_{B_1(x)}(E)$  of a set $E$ relative to $B_1(x)$ is defined as
\begin{equation}\label{cap}
{\rm cap}_{B_1(x)}(E) = \inf \bigg\{\int _{B_1(x)} |\nabla v|^2\, dy : v\in C^{0,1}_0(B_1(x)),\, v\geq 1 \,\hbox{on}\, E\bigg\}\,,
\end{equation}
where  $C^{0,1}_0(B_1(x))$ denotes the space of  Lipschitz continuous, compactly supported functions in $B_1(x)$. 
Let us notice that, 
if $n \geq 3$, then the  capacity ${\rm cap}_{B_1(x)}$  is equivalent (up to multiplicative constants depending on $n$) to the standard capacity in the whole of $\rn$.

A more transparent condition on $\partial \o$, though slightly stronger  than  \eqref{capcond}, for our regularity result to hold can be given in terms of an integrability property of $\mathcal B$. It involves  the space of weak type, also called Marcinkiewicz space, with respect to the measure $\hh$ on $\partial \o$ defined as
\begin{equation}\label{Xspace}
X= \begin{cases} L^{n-1, \infty} & \quad \hbox{if $n \ge 3$,}
\\
L^{1, \infty} \log L  & \quad \hbox{if $n =2$.}
\end{cases}
\end{equation}
Here,  $L^{n-1, \infty}$ denotes the
weak Lebesgue space endowed with the norm
\begin{equation}\label{weakleb}
\|h\|_{L^{q, \infty} (\partial \Omega)} = \sup _{s \in (0, \hh(\partial \Omega))} s ^{\frac 1q} h^{**}(s),
\end{equation}
and $L^{1, \infty}\log L$ the weak Zygmund space endowed with the norm
\begin{equation}\label{weaklog}
\|h\|_{L^{1, \infty} \log L (\partial \Omega)} = \sup _{s \in (0, \hh(\partial \Omega))} s \log\big(1+ \tfrac{C}s\big) h^{**}(s),
\end{equation}
for any constant  $C>\hh(\partial \Omega)$. Observe that  different constants $C$ result in equivalent norms in \eqref{weaklog}.
 Here, $h ^{**}(s)= \smallint _0^s h^* (r)\, dr$ for $s >0$, where $h^*$ denotes the decreasing rearrangement of a measurable function $h: \partial \o \to \mathbb R$ with respect to $\hh$.
%
\\ The relevant integrability property on $\mathcal B$ amounts to requiring that $\mathcal B \in X(\partial \o)$, an assumption that will be denoted by $\partial \o \in W^2X$.  An additional  smallness condition is however again needed on local norms of $\mathcal B$ in $X(\partial \o)$, which takes the form
\begin{equation}\label{smalln}
 \lim _{r\to 0^+} \Big(\sup _{x \in \partial \Omega} \|\mathcal B \|_{X(\partial \Omega \cap B_r(x))}\Big) < c_2 \,,
\end{equation}
for a suitable constant $c_2=c_2(n, N, p, L_\Omega, d_\Omega)$.
\par Observe that condition \eqref{smalln} is certainly  fulfilled if either $n\geq 3$ and $\partial \o \in W^{2,{n-1}}$, or $n=2$ and  $\partial \o \in W^{2, q}$ for some $q>1$,  and hence, in particular, if $\partial \o \in C^2$. Let us however emphasize that 
 condition \eqref{smalln} does not even entail that $\partial \o \in C^1$. 
 \par The link bewteen assumptions \eqref{smalln} and \eqref{capcond} is provided by
 \cite[Lemmas 3.5 and 3.7]{cmsecondsys}.  Those results  ensure that, given any constant $c_1$, there exists a constant $c_2=c_2(n, d_\o, L_\o, c_1)$ such that if $\o$ fulfills  \eqref{smalln}, then it also satisfies \eqref{capcond}.

\begin{theorem}\label{secondcap} {\rm {\bf [Second-order estimates under minimal boundary regularity]}}
Let $N=1$, $1<p <\infty$ and  $T>0$. 
Assume that
 $\o$ is a bounded  Lipschitz  domain in $\rn$, $n \geq 2$, with $\partial \o \in W^{2,1}$. 
Let $f \in
L^2(\Omega_T )$ and $\psi \in W^{1,p}_0(\o)$, and let $u$ be the approximable solution to problem \eqref{eqdirichlet}.
There exists a constant $c_1=c_1(n,p, d_\o, L_\o)$ such that, if  condition \eqref{capcond} is fulfilled, then
\begin{equation}\label{secondcap1} |\nabla u|^{p-2} \nabla u
\in L^2((0,T); W^{1,2}(\o))\,.
\end{equation}
Moreover, there exists a constant $C=C(n,p, \o, T)$ such that
\begin{align}\label{secondcap2} 
 \||\nabla u|^{p-2} \nabla u\|_{L^2((0,T); W^{1,2}(\o))} \leq C \big(\|f\|_{L^2(\o_T)} + \|\nabla \psi\|_{L^p(\o)}^{\frac p2}\big)\,.
\end{align}
%
%
%
%
In particular, there exists a constant  $c_2=c_2(n,p, d_\o, L_\o)$ such that properties \eqref{secondcap1} and \eqref{secondcap2} hold if $\partial \o \in W^2X$ and fulfills condition \eqref{smalln}.
\end{theorem}

The sharpness of Theorem \ref{secondcap}, both in the  assumptions on $\o$ and in the estimate for $u$, is pointed out in the following remarks.

\begin{remark}\label{twosided}  {\rm {\bf [Sharpness of estimates]}}
{\rm 
The  regularity for the derivative in $t$ of $u$ and the derivatives in $x$ of $|\nabla u|^{p-2} \nabla u$, provided by Theorems \ref{exist} and \ref{secondcap}, is sharp. Indeed, assume, for instance, that $\psi =0$. One trivially has that
\begin{align}\label{twosided1}
\|f\|_{L^2(\o_T)}  & = \|u_t - {\rm div} (|\nabla u|^{p-2} \nabla u )  \|_{L^2(\o_T)}  
 \leq \|u_t \|_{L^2(\o_T)} +  \|{\rm div} (|\nabla u|^{p-2} \nabla u )  \|_{L^2(\o_T)} 
\\ \nonumber & \leq 
 \|u_t \|_{L^2(\o_T)} + c \||\nabla u|^{p-2} \nabla u\|_{L^2((0,T); W^{1,2}(\o))} 
\end{align}
for some constant $c=c(n,p)$. Thus,  if $\o$ is as in Theorem \ref{secondcap}, then the following two-sided estimate holds:
\begin{align}\label{twosided2}
c_1 \|f\|_{L^2(\o_T)} \leq    \|u_t \|_{L^2(\o_T)} +  \||\nabla u|^{p-2} \nabla u\|_{L^2((0,T); W^{1,2}(\o))} 
  \leq c_2  \|f\|_{L^2(\o_T)}
\end{align}
for suitable positive constants $c_1$ and $c_2$.
}
\end{remark}

\begin{remark}\label{sharp}  {\rm {\bf [Sharpness of boundary regularity]}}
{\rm 
Membership of $\partial \o$ in $W^2X$, namely the mere finiteness of  the limit in \eqref{capcond}, is not sufficient for the conclusions of Theorem \ref{secondcap} to hold.
Actually,  in \cite{Ma67, Ma73} a family of domains $\{\o_\beta\}$ with the following properties is exhibited: (i) the limit in \eqref{capcond}, with $\o = \o_\beta$, is finite and depends continuously on $\beta$; (ii)  
the stationary solution in $\Omega _\beta$ to the heat equation,
with a suitable smooth right-hand side,   belongs to $W^{2,2}(\o_T)$ if and only if the relevant limit  does not exceed an explicit threshold.
%
%
The boundary of each domain $\o_\beta$ is smooth outside a small region, where it agrees with the graph of a function $\Theta_\beta$ depending on the variables $(x_1, \dots , x_{n-1})$ only through $x_1$ and having the form 
\begin{equation}\label{theta}
\Theta_\beta (x_1, \dots , x_{n-1})= \beta |x_1|(\log |x_1|)^{-1}
\end{equation}
  for small $x_1$. The limit  in  \eqref{capcond}   is a multiple of  $\beta$, and the  solution $u$ turns out to belong to $W^{2,2}(\o)$ if and only if the constant $\beta$ is smaller than or equal to an explicit value, depending only on $n$.
\\
The sharpness of the alternate
assumption \eqref{smalln} in Theorem \ref{secondcap}   can be shown, for instance, when $n=3$ and   $p\in (\tfrac 32,2]$, by an 
 example from \cite{KrolM}, which again applies to the stationary case.  In that paper,  open sets  $\o \subset \mathbb R^3$ are constructed such that  the limit in \eqref{smalln} is finite,  but too large,  and  the stationary solution $u$ to problem \eqref{eqdirichlet}, with  a smooth right-hand side,  is so irregular that $|\nabla u|^{p-2}\nabla u \notin W^{1,2}(\Omega_T)$.   Similarly, in $\mathbb R^2$ there exist open sets $\Omega$
for which the limit in  \eqref{smalln}   is finite but larger than some critical value,
and  where the stationary solution $u$ to  the heat equation  with  a smooth right-hand side fails to belong to $W^{2,2}(\o)$ \cite{Ma67}. 
\\ Let us notice that, if $\o_\beta$ is a domain as above, then  $\partial \o_\beta \notin   W^2L^{2, \infty}$ if $n \geq 3$.  Hence, if $\beta $ is sufficiently small, the capacitary criterion \eqref{capcond}  of Theorem \ref{secondcap} applies to deduce properties \eqref{secondcap1} and \eqref{secondcap2}, whereas the integrability condition \eqref{smalln} does not.
 }
\end{remark}

\begin{remark}\label{1-dim}  {\rm {\bf [Case $n=1$]}}
{\rm An inspection of the proof will reveal that, if $n=1$, then Theorem \ref{secondcap} holds  if $\o$ is any bounded interval. The result is considerarbly easier in this case, since the divergence operator agrees with plain differentiation in dimension one. }
\end{remark}

Under the assumption that $\o$ is  a bounded convex open set, the conclusions of Theorem \ref{secondcap} hold without any additional regularity condition on $\partial \o$.  

\begin{theorem}\label{secondconvex} {\rm {\bf [Second-order estimates in convex domains]}}
Let $N=1$, $1 <p <\infty$ and $T>0$. Assume that $\o$ is a bounded  convex open set in $\rn$, $n \geq 2$. Let $f \in
L^2(\Omega_T )$ and $\psi \in W^{1,p}_0(\o)$, and let $u$ be the approximable solution to problem \eqref{eqdirichlet}. Then
$ |\nabla u|^{p-2} \nabla u
\in L^2((0,T); W^{1,2}(\o))$ and inequality \eqref{secondcap2}  holds for some positive constant  $C=C(n,p, \o, T)$.
\end{theorem}


Our results for systems are stated in Theorems \ref{secondcapsys} and \ref{secondconvexsys} below. They parallel Theorems  \ref{secondcap} and \ref{secondconvex}. The difference here is that they hold for $p \geq 2$ and can hence be stated in terms of weak solutions to problem \eqref{eqdirichlet}. Of course, Remarks \ref{twosided}--\ref{1-dim} carry over to this case.

\begin{theorem}\label{secondcapsys} {\rm {\bf [Second-order estimates for systems under minimal boundary regularity]}}
Let $N\geq 1$, $2 \leq p < \infty$ and  $T>0$. 
Assume that
 $\o$ is a bounded  Lipschitz  domain with $\partial \o \in W^{2,1}$. 
Let $f \in
L^2(\Omega_T )$ and $\psi \in W^{1,p}_0(\o)$, and let $u$ be the weak solution to problem \eqref{eqdirichlet}.
There exists a constant $c_1=c_1(n, N, p, d_\o, L_\o)$ such that, if  condition \eqref{capcond} is fulfilled, then
\begin{equation}\label{secondcapsys1} |\nabla u|^{p-2} \nabla u
\in L^2((0,T); W^{1,2}(\o))\,.
\end{equation}
Moreover, there exists a constant $C=C(n,N, p, \o, T)$ such that
\begin{align}\label{secondcapsys2} 
 \||\nabla u|^{p-2} \nabla u\|_{L^2((0,T); W^{1,2}(\o))} \leq C \big(\|f\|_{L^2(\o_T)} + \|\nabla \psi\|_{L^p(\o)}^{\frac p2}\big)\,.
\end{align}
In particular, there exists a constant  $c_2=c_2(n, N, p, d_\o, L_\o)$ such that properties \eqref{secondcapsys1} and \eqref{secondcapsys2} hold if $\o \in W^2X$ and fulfils condition \eqref{smalln}.
\end{theorem}

\begin{theorem}\label{secondconvexsys} {\rm {\bf [Second-order estimates for systems in convex domains]}}
Let $N\geq 1$, $2\leq p <\infty$ and $T>0$. Assume that $\o$ is a bounded  convex open set in $\rn$, $n \geq 2$. Let $f \in
L^2(\Omega_T)$ and $\psi \in W^{1,p}_0(\o)$, and let $u$ be the weak solution to problem \eqref{eqdirichlet}. Then
$ |\nabla u|^{p-2} \nabla u
\in L^2((0,T); W^{1,2}(\o))$ and inequality \eqref{secondcapsys2}  holds for some positive constant  $C=C(n, N, p, \o, T)$.
\end{theorem}

\section{Proofs}\label{proofs}

Our proof of Theorem \ref{exist} combines arguments introduced in \cite{BBGGPV} in the elliptic case and developed in \cite{Pr} for parabolic problems, with an estimate for $u_t$ and $\nabla u$ in the spirit of  \cite[Proposition 4.1]{FS}, and of earlier results from \cite[Chapter 1, Theorem 8.1]{lions} and \cite{BEKP}.

\medskip
\par\noindent
{\bf Proof  of Theorem \ref{exist}}.  We assume that $n \geq 2$, the case when $n=1$ being analogous, and even simpler.
Let $\{f_k\} \subset C^\infty_0(\o_T)$ and $\psi_k \in C^{\infty}_0(\o)$ be sequences such that
\begin{equation}\label{aug1}
f_k \to f \quad \hbox{in $L^2(\o_T)$} \quad \hbox{and} \quad \|f_k\|_{L^2(\o_T)} \leq 2  \|f\|_{L^2(\o_T)}\,,
\end{equation}
\begin{equation}\label{aug1bis}
\psi_k \to \psi \quad \hbox{in $W^{1,p}(\o)$} \quad \hbox{and} \quad \|\nabla \psi_k\|_{L^{p}(\o)} \leq 2  \|\nabla \psi\|_{L^{p}(\o)}\,,
\end{equation}
and let $\{u_k\}$ be the corresponding sequence of weak solutions to problems \eqref{eqdirichletk}. 
A global in time version of \cite[Proposition 4.1]{FS}  tells us that  $u_k \in L^\infty((0, T); W^{1,p}(\o))$, $(u_k)_t \in L^2(\o_T)$ and 
\begin{align}\label{solamaink} 
\| u_k\|_{L^\infty((0,T); W^{1,p}_0(\o))}^{\frac p2}+ \|(u_k)_t\|_{L^2(\o_T)}& \leq C \big(\|f_k\|_{L^2(\o_T)} + \|\nabla \psi_k\|_{L^p(\o)}^{\frac p2}\big) \\ \nonumber & \leq C'  \big(\|f\|_{L^2(\o_T)} + \|\nabla \psi\|_{L^p(\o)}^{\frac p2}\big)  \,,
\end{align}
for suitable constants $C$ and $C'$ depending on $n$, $p$, $|\o|$ and $T$.
Inequality \eqref{solamaink} holds in the whole interval $[0,T]$, namely up to $t=0$, instead of just locally, as in the result of \cite{FS},  thanks to the present assumption  that $\psi$ belongs  to $W^{1,p}_0(\o)$.  The proof  is completely analogous to (in fact, simpler than) that of \cite[Proposition 4.1]{FS},  and will be omitted. 
\\ One clearly has that $L^\infty((0,T); W^{1,p}_0(\o)) \to L^p(\o_T)$. 
Thus, owing to \eqref{solamaink}, the sequence $\{u_k\}$ is bounded in the anisotropic Sobolev space $W^{1,(p,2)}(\o_T)$ defined as
$$W^{1,(p,2)}(\o_T)=\{v:\, v\in L^p(\o_T), \,|\nabla v|\in L^p(\o_T), v_t \in L^2(\o_T)\}\,,$$
and equipped with the norm 
$$\|v\|_{W^{1,(p,2)}(\o_T)} = \|v\|_{L^p(\o_T)} + \|\nabla v\|_{L^p(\o_T)}+ \|v_t\|_{L^2(\o_T)}\,.$$
Since $p \in (1,\infty)$, the space $W^{1,(p,2)}(\o_T)$ is reflexive. Moreover, on setting
$q=\min\{p, 2\}$, one trivially has $W^{1,(p,2)}(\o_T)\to W^{1,q}(\o_T)$, and hence the sequence $\{u_k\}$ is also bounded in $W^{1,q}(\o_T)$. 
Therefore, given any number  $r \in \big[1, \tfrac{q(n+1)}{n+1-q}\big)$,
there exists 
%
%
a subsequence, still denoted by $\{u_k\}$, such that
\begin{equation}\label{aug2}
u_k \rightharpoonup u \quad \hbox{weakly in $W^{1, (p,2)}(\o_T)$} \quad \hbox{and} \quad u\to u \quad \hbox{in $L^r(\o_T)$.}
\end{equation}
\\ Our goal is now to show that $\{\nabla u_k\}$ is a Cauchy sequence in measure. Fix any $\lambda >0$ and $\ep>0$. Given any $\theta , \sigma >0$, one has that
\begin{align}\label{aug3}
|\{|\nabla u_k -\nabla u_m|>\lambda\}| & \leq |\{|\nabla u_k|>\theta\}| +|\{|\nabla u_m|>\theta\}| + |\{|u_k -u_m|>\sigma\}|\\ \nonumber & \quad + |\{|\nabla u_k -\nabla u_m|>\lambda, |u_k -u_m|\leq\sigma, |\nabla u_k|\leq \theta, |\nabla u_m|\leq\theta \}|
\end{align}
for $k,m \in \mathbb N$. Since the sequence $\{u_k\}$ is bounded in $W^{1, (p,2)}(\o_T)$, there exists $\theta >0$ such that
\begin{align}\label{aug4}
|\{|\nabla u_k|>\theta\}| +|\{|\nabla u_m|>\theta\}|< \ep
\end{align}
for every $k, m \in \mathbb N$. Define, for $\sigma >0$, the function $T_\sigma : \mathbb R \to \mathbb R$ as
$$T_\sigma (s) = \begin{cases} s & \quad \hbox{if $|s|< \sigma$,}
\\ \sigma \,{\rm sign} (s) & \quad \hbox{if $|s|\geq \sigma$.}
\end{cases}
$$
Clearly, $|T_\sigma (s)|\leq \sigma $ for $s \in \mathbb R$.
 Choose the test function $\phi = T_\sigma (u_k-u_m)$ in the weak formulation of problem \eqref{eqdirichletk}, namely in equation \eqref{lions} with $u$, $f$, $\psi$ replaced by $u_k$, $f_k$, $\psi_k$. Next, choose the same test function in the same problem with $k$ replaced by $m$, and subtract the resultant equations to obtain that
\begin{align}\label{aug5}
\int _0^T   & \langle (u_k - u_m)_t ,T_\sigma (u_k - u_m) \rangle \,dt \\ \nonumber & \quad +  \int _0^T \int_\o (|\nabla u_k|^{p-2}\nabla u_k - |\nabla u_m|^{p-2}\nabla u_m)\cdot \nabla (T_\sigma (u_k - u_m))\, dx\,dt \\ \nonumber & 
= \int _0^T \int_\o  (f_k - f_m)T_\sigma (u_k - u_m)\, dx\,dt \,.
\end{align}
 Define 
the function $\Lambda_\sigma : \mathbb R \to \mathbb R$ as $\Lambda _\sigma (s) = \int _0^sT_\sigma (r)\, dr$ for $s \in \mathbb R$. Observe that 
\begin{equation}\label{sep1}
0\leq \Lambda _\sigma (s) \leq \sigma |s|  \quad \hbox{for $s\in \mathbb R$.}
\end{equation}
 Moreover, one has that 
\begin{equation}\label{aug30}
\int _0^\tau \langle v_t , T_\sigma (v) \rangle \, dt  = \int_\o \Lambda _\sigma  (v(x,\tau)\, dx -   \int_\o \Lambda _\sigma  (v(x,0)\, dx\quad \hbox{if $\tau  \in(0,T)$,}
\end{equation}
provided that $v\in L^2(\o_T)$ and $v_t  \in L^2(\o_T)$ -- see e.g. by \cite[Section 2.1.1]{GaMa}. 
From equations \eqref{aug5}, \eqref{sep1}    and  \eqref{aug30}, with $v=u_k -u_m$,  one can deduce that
\begin{align}\label{aug6}
\int _0^T &\int_{\{|u_k-u_m|<\sigma\}}(|\nabla u_k|^{p-2}\nabla u_k - |\nabla u_m|^{p-2}\nabla u_m)\cdot (\nabla u_k - \nabla u_m)\, dxdt \\ \nonumber & \leq   
 \sigma \int _0^T \int_\o  |f_k - f_m| \, dxdt + \sigma \int _\o |\psi _k - \psi_m|\, dx \leq \sigma C (\|f\|_{L^2(\o_T)} + \|\nabla \psi\|_{L^p(\o)}) \,
\end{align}
for  some constant $C=C(n,p, |\o|)$ and every $k, m \in \mathbb N$. Define
$$\kappa = \min \{(|\xi|^{p-2}\xi - |\eta|^{p-2}\eta) \cdot (\xi -\eta): |\xi|\leq \theta, |\eta|\leq \theta, |\xi-\eta|\geq \lambda\}\,.$$
Since $\kappa >0$, from inequality \eqref{aug6} one infers that
\begin{multline}\label{aug7}
\kappa \, |\{|\nabla u_k -\nabla u_m|>\lambda, |u_k -u_m|\leq\sigma, |\nabla u_k|\leq \theta, |\nabla u_m|\leq\theta \}| \\ \leq \sigma C (\|f\|_{L^2(\o_T)} + \|\nabla \psi\|_{L^p(\o)})\,,
\end{multline}
whence
\begin{align}\label{aug8}
|\{|\nabla u_k -\nabla u_m|>\lambda, |u_k -u_m|\leq\sigma, |\nabla u_k|\leq \theta, |\nabla u_m|\leq\theta \}| < \ep
\end{align}
for every $k, m \in \mathbb N$, provided that $\sigma $ is sufficiently small. On the other hand, by property \eqref{aug2}, $\{u_k\}$ is a Cauchy sequence in measure, and hence
\begin{align}\label{aug9}
|\{|u_k-u_m|>\sigma\}|< \ep
\end{align}
if $k$ and $m$ are sufficiently large. Inequalities \eqref{aug3}, \eqref{aug4}, \eqref{aug8} and \eqref{aug9} tell us that $\{\nabla u_k\}$ is actually a Cauchy sequence in measure. As a consequence, there exists a subsequence, still indexed by $k$, such that
\begin{equation}\label{aug10}
\nabla u_k \to \nabla u \quad \hbox{a.e. in $\o_T$.}
\end{equation}
Hence, equation \eqref{june0} holds. Moreover, passing to the limit  in inequality \eqref{solamaink} as $k \to \infty$, yields inequality \eqref{solamain}.
\\ In order to establish property \eqref{june0bis}, fix any $\tau \in [0,T]$ and make use of the test function $\phi= T_1(u_k-u_m)\chi _{[0, \tau)}$ in the weak formulation of problem \eqref{eqdirichletk}  and in its analogue with $k$ replaced by $m$. Here $\chi_E$ stands for the characteristic function of the set $E$.  Subtracting the resultant equations yields 
\begin{align}\label{aug11}
\int _0^\tau & \langle (u_k - u_m)_t , T_1(u_k - u_m)\rangle \, dt\\  \nonumber & \quad  +  \int _0^\tau \int_\o (|\nabla u_k|^{p-2}\nabla u_k - |\nabla u_m|^{p-2}\nabla u_m)\cdot \nabla (T_1 (u_k - u_m))\, dx\,dt \\ \nonumber & 
= \int _0^\tau \int_\o  (f_k - f_m)T_1 (u_k - u_m)\, dx\,dt \,.
\end{align}
Since the integrand in the second integral on the left-hand side of inequality \eqref{aug11} is nonnegative, one  infers via equation \eqref{aug30}, with $v= u_k -u_m$,  and \eqref{sep1} that
\begin{align}\label{aug12}
\int_\o \Lambda_1(u_k(x,\tau) - u_m(x,\tau))\, dx & \leq  \int_\o \Lambda_1(\psi_k  - \psi_m )\, dx  + 
\int _0^\tau \int_\o  (f_k - f_m)T_1 (u_k - u_m)\, dxdt
\\ \nonumber & \leq \|\psi_k - \psi _m\|_{L^1(\o)} +
T  \|f_k - f_m\|_{L^1(\o_T)}  \\ \nonumber &\leq C( \|\nabla \psi_k - \nabla \psi _m\|_{L^p(\o)} +
 \|f_k - f_m\|_{L^2(\o_T)}) \,,
\end{align}
for some constant $C=C(n,p,|\o|, T)$ and every $\tau \in [0,T]$. On the other hand, 
\begin{align}\label{aug13}
& \int_\o  \Lambda_1 (u_k(x,\tau) - u_m(x,\tau))\, dx  \\  \nonumber & \geq \int _{\{|u_k-u_m|\leq 1\}} |u_k(x,\tau)-u_m(x,\tau)|^2\, dx + \frac 12\int _{\{|u_k-u_m|> 1\}} |u_k(x,\tau)-u_m(x,\tau)|\, dx
\\ \nonumber & \geq 
\frac 1{|\o|}\bigg(\int _{\{|u_k-u_m|\leq 1\}} |u_k(x,\tau)-u_m(x,\tau)|\, dx\bigg)^2 + \frac 12\int _{\{|u_k-u_m|> 1\}} |u_k(x,\tau)-u_m(x,\tau)|\, dx
\end{align}
for every $\tau \in [0,T]$.
Inequalities \eqref{aug12} and \eqref{aug13} ensure that $u_k$ is a Cauchy sequence in $C([0,T]; L^1(\o))$, and hence that  \eqref{june0bis} holds (up to subsequences).
\\ Finally, as far as uniqueness is concerned, assume that $u$ and $\overline u$ are   approximable solutions to problem \eqref{eqdirichlet}.
Let $\{f_k\}$, $\{\psi_k\}$ and $\{u_k\}$ be sequences as in the definition of approximable solution for $u$, and $\{\overline f_k\}$, $\{\overline \psi_k\}$ and $\{\overline u_k\}$ sequences as in a parallel definition for $\overline u$. On choosing, for $\sigma >0$, the test function $\phi = T_\sigma (u_k - \overline u_k)$ in the definitions of weak solutions for $u_k$ and $\overline u_k$,  and subtracting the equations so obtained, one deduces, analogously to \eqref{aug6},
\begin{align}\label{aug14}
\int _0^T &\int_{\{|u_k-\overline u_k|<\sigma\}}(|\nabla u_k|^{p-2}\nabla u_k - |\nabla \overline u_k|^{p-2}\nabla \overline u_k)\cdot (\nabla u_k - \nabla \overline u_k)\, dx\,dt \\ \nonumber & \leq   
 \sigma \int _0^T \int_\o  |f_k - \overline f_k| \, dxdt + \sigma \int _\o |\psi _k - \overline \psi_k|\, dx
\end{align}
for $k \in \mathbb N$. By our assumptions, the right-hand side of inequality \eqref{aug14} approaches $0$ as $k \to \infty$, and $u_k \to u$, $\nabla u_k \to \nabla u$, $\overline u_k \to \overline u$ and $\nabla \overline u_k \to \overline u$ a.e. in $\o_T$. Hence, inequality  \eqref{aug14} implies, via Fatou's lemma, that
\begin{align}\label{aug15}
\int _0^T &\int_{\{|u-\overline u|<\sigma\}}(|\nabla u|^{p-2}\nabla u - |\nabla \overline u|^{p-2}\nabla \overline u)\cdot (\nabla u - \nabla \overline u)\, dx\,dt =0
\end{align}
for every $\sigma >0$. The integrand in \eqref{aug15} is nonnegative, and vanishes if and only if $\nabla u = \nabla \overline u$. Hence, by the arbitrariness of $\sigma$, we have that $\nabla u = \nabla \overline u$ a.e. in $\o_T$. Inasmuch as $u$ and $\overline u \in L^1((0,T); W^{1,1}_0(\o))$, we conclude that $u=\overline u$. \qed

We are now ready to prove Theorem \ref{secondcap}.

\medskip
\par\noindent
{\bf Proof of Theorem \ref{secondcap}}.
Assume, for the time being, 
that 
\begin{equation}\label{smooth}
f \in C^\infty _0 (\Omega_T) \quad \hbox{and} \quad \psi \in C^\infty _0(\Omega)\,.
\end{equation}
Then, there exists a unique weak solution $u \in C([0,T]; L^2(\o))\cap L^p((0,T); W^{1,p}_0(\o))$, with $u_t \in L^{p'}((0,T); W^{-1,p'}(\o))$,  to problem \eqref{eqdirichlet} in the sense of \eqref{lions}.
By inequality \eqref{solamaink}, with $f_k$, $\psi_k$ and $u_k$ replaced with  $f$, $\psi$ and $u$, we also have that $u_t \in L^2(\o _T)$.
  Choose any test function $\phi$ of the form $\phi (x, t) = \varphi (x)\rho(t)$ in \eqref{lions}, with $\varphi \in W^{1,p}_0(\o)\cap L^\infty(\o)$ and $\rho \in C^\infty_0(0,T)$. Equation 
\eqref{lions}
then entails that
%
\begin{equation}\label{june2}
 \int _0^T \rho(t)\int _\o |\nabla u(x,t)|^{p-2}\nabla u(x,t) \cdot \nabla \varphi (x)\, dx \,dt = \int _0^T\rho (t)\int_\o \varphi (x) (f(x,t) - u_t(x,t))\, dx \, dt\,.
\end{equation}
Hence, by  the arbtrariness of the function $\rho$, 
\begin{equation}\label{june3}
\int _\o |\nabla u(x,t)|^{p-2}\nabla u(x,t) \cdot \nabla \varphi (x)\, dx = \int_\o \varphi (x) (f(x,t) - u_t(x,t))\, dx \quad \hbox{for a.e. $t \in (0,T)$,}
\end{equation}
and for every $\varphi \in W^{1,p}_0(\o)\cap L^\infty(\o)$. Now, define, 
 for each $t \in (0,T)$, the function $g^t: \o \to \mathbb R$ as 
\begin{equation}\label{g}
g^t(x) =f(x,t) - u_t(x,t) \quad \hbox{for $x \in \o$,}
\end{equation}
and the function $w^t : \o \to \mathbb R$ as 
\begin{equation}\label{w}
w^t(x) = u(x,t) \quad \hbox{for $x \in \o$.}
\end{equation}
Then $g^t \in L^2(\o)$ and $w^t \in W^{1,p}_0(\o)$ for a.e. $t \in (0, T)$. Moreover, equation \eqref{june3} reads
\begin{equation}\label{june3bis}
\int _\o |\nabla w^t|^{p-2}\nabla w^t \cdot \nabla \varphi \, dx = \int_\o \varphi \,  g^t \, dx
\end{equation}
for every $\varphi \in W^{1,p}_0(\o)\cap L^\infty(\o)$ and for a.e. $t \in (0, T)$. Fix any such $t$, and consider the (elliptic) Dirichlet problem 
\begin{equation}\label{aux}
\begin{cases}
- {\rm div} (|\nabla v^t|^{p-2}\nabla v^t ) = g^t  & {\rm in}\,\,\, \o \\
 v^t =0  &
{\rm on}\,\,\,
\partial \o \,.
\end{cases}
\end{equation}
This problem has a generalized approximable solution $v^t$ in the following sense. For every $\sigma >0$, the function $T_\sigma (v^t) \in W^{1,p}_0(\o)$,  and there exists a measurable function $Z^t : \o \to \rn$ such that
\begin{equation}\label{sep2} \nabla T_\sigma (v^t)  = Z^t \chi_{\{|v^t|<\sigma\}} \quad \hbox{a.e. in $\o$,}
 \end{equation}
 for every $\sigma >0$. Moreover, there exist sequences
$\{g_k^t\} \subset C^\infty_0(\o)$ and $\{v_k^t\} \subset W^{1,p}_0(\o)$ such that $g_k^t \to g$ in $L^1(\o)$, 
 $v_k^t$ is the weak solution to problem 
\begin{equation}\label{L1k}
\begin{cases}
- {\rm div} (|\nabla v_k^t|^{p-2} \nabla v_k^t ) = g_k^t  & {\rm in}\,\,\, \o \\
 v_k^t =0  &
{\rm on}\,\,\,
\partial \o\,,
\end{cases}
\end{equation}
and 
\begin{equation}\label{L1conv}
v_k^t \to v^t \quad \hbox{and} \quad \nabla v_k \to Z^t \quad \hbox{a.e. in $\o$.}
\end{equation}
This follows e.g. from   \cite[Theorem 3.2]{cmapprox}. We claim that
\begin{equation}\label{june4}
w^t= v^t\,.
\end{equation}
In order to verify this claim, observe that $T_\sigma(w^t-v_k^t) \in W^{1,p}_0(\Omega) \cap L^\infty (\o)$ for every $k\in \mathbb N$ and $\sigma>0$,
since   $w^t \in W^{1,p}_0(\Omega)$ and $v_k^t \in W^{1,p}_0(\Omega)$. Making use of the test function $\varphi = T_\sigma(w^t-v_k^t)$   in equation \eqref{june3bis} and in the definition of weak solution to problem \eqref{L1k}, and subtracting the resultant equations tells us that
\begin{equation}\label{sola3}
\int_{\{|w^t-v_k^t|< \sigma\}}(|\nabla w^t|^{p-2}\nabla w^t - |\nabla v_k^t|^{p-2}\nabla v_k^t)\cdot (\nabla w^t - \nabla v_k^t)\, dx = \int _\o (g^t -g_k^t)\, T_\sigma(w^t-v_k^t) \, dx\,.
\end{equation}
Since $|T_\sigma(w^t-v_k^t)|\leq \sigma$ and $g_k^t \to g^t$ in $L^1(\o)$, the right-hand side of equation \eqref{sola3} tends to $0$ as $k \to \infty$. Thus, inasmuch as  the integrand on the left-hand side  is nonnegative,  passing to the limit in \eqref{sola3} as $k \to \infty$ yields, by Fatou's lemma and equation \eqref{L1conv},
\begin{equation}\label{sola4}
\int_{\{|w^t-v^t|< \sigma\}}(|\nabla w^t|^{p-2}\nabla w^t - |Z^t|^{p-2}Z^t)\cdot (\nabla w^t - Z^t)\, dx =0
\end{equation}
for every $\sigma>0$. The integrand in \eqref{sola4} in nonnegative, and vanishes if and only if $\nabla w^t = Z^t$. Therefore,  for every $\sigma>0$,  we have that $\nabla w^t = Z^t$ a.e. in the set $\{|w^t-v^t|< \sigma\}$, whence 
\begin{equation}\label{sola5}
\nabla w^t = Z^t \qquad \hbox{a.e. in $\o$.}
\end{equation}
Now, the function $T_\sigma (w^t-T_\varrho(v^t)) \in W^{1,p}_0(\o)\subset W^{1,1}_0(\o)$ for every $\sigma, \varrho >0$. An application of the Sobolev inequality in $W^{1,1}_0(\o)$ and the use of equations \eqref{sep2} and  \eqref{sola5} imply that
\begin{equation}\label{sola6}
\bigg( \int_{\o }|T_{\sigma}(w^t-T_\varrho(v^t)|^{n'}\, dx \bigg
)^{n'}
\le C
\bigg(
\int_{\{\varrho<|w^t|<\varrho+\sigma\} }|\nabla w^t|\, dx+
\int_{\{\varrho-\sigma<|w^t|<\varrho\} }|\nabla w^t|\, dx
\bigg)
\end{equation}
for some constant $C=C(n)$ and for every $\sigma, \varrho >0$. Since $|\nabla w^t| \in L^1(\o)$, for each fixed $\sigma >0$ the right-hand side of inequality \eqref{sola6} converges to $0$ as $\varrho\to \infty$. Therefore,  passing to the limit  in \eqref{sola6} as $\varrho\to \infty$ tells us that
$$\int_{\o }|T_{\sigma}(w^t-v^t)|^{n'}\, dx =0$$
for every $\sigma >0$. Hence, equation \eqref{june4} follows, on letting $\sigma \to \infty$.
An application of \cite[Theorem 2.1]{cmsecondsys} tells us that
\begin{equation}\label{aug25}
|\nabla v^t|^{p-2} \nabla v^t \in W^{1,2}(\o),\end{equation}
 and there exists a constant $C=C(n,p, \o)$ such that
\begin{equation}\label{aug26}
\||\nabla v^t|^{p-2} \nabla v^t\|_{W^{1,2}(\o)}^2 \leq C \|g^t\|_{L^2(\o)}.
\end{equation}
Hence, owing to equations \eqref{june4}, \eqref{g} and \eqref{w}  
$$ |\nabla u(\cdot , t)|^{p-2} \nabla u(\cdot , t) \in W^{1,2}(\o)$$ 
for a.e. $t \in (0,T)$,
and 
\begin{equation}\label{june5}
\||\nabla u(\cdot , t)|^{p-2} \nabla u(\cdot , t)\|_{W^{1,2}(\o)}^2 \leq C \int_\o f(x,t)^2 + u_t(x,t)^2\, dx 
\end{equation}
for the same constant  $C=C(n,p, \o)$.
Integrating inequality \eqref{june5} with respect to $t$ over $(0,T)$ yields
\begin{equation}\label{june6}
\||\nabla u|^{p-2} \nabla u\|_{L^2((0,T); W^{1,2}(\o))}^2 \leq C \big(\|f\|_{L^2(\o_T)}^2 + \|u_t\|_{L^2(\o_T)}^2\big)\,.
\end{equation}
Coupling inequalities  \eqref{june6} and \eqref{solamain} tells us that $|\nabla u|^{p-2} \nabla u \in L^2((0,T); W^{1,2}(\o))$, and 
\begin{align}\label{june7}
\||\nabla u|^{p-2} \nabla u\|_{L^2((0,T); W^{1,2}(\o))} \leq C \big(\|f\|_{L^2(\o_T)} + \|\nabla \psi\|_{L^p(\o)}^{\frac p2}\big)\,
\end{align}
constant  $C=C(n,p, \o)$.
\\
We have now to remove the additional assumptions  \eqref{smooth}. To this purpose, consider sequences  $\{f_k\}$, $\{\psi_k\}$ and $\{u_k\}$  as in the definition of approximable solution to problem \eqref{eqdirichlet}. Clearly, we may assume that $\|f_k\|_{L^2(\o_T)} \leq 2 \|f\|_{L^2(\o_T)}$ and $\|\nabla \psi _k\|_{L^p(\o)} \leq 2\|\nabla \psi \|_{L^p(\o)}$ for $k \in \mathbb N$.
By inequality    \eqref{june7}, applied with $u$ replaced by $u_k$, and our assumptions on the sequences $\{f_k\}$ and $\{\psi _k\}$, there exists a constant $C=C(n,p, \o)$  such that
\begin{align}\label{june8}
 \||\nabla u_k|^{p-2} \nabla u_k\|_{L^2((0,T); W^{1,2}(\o))}
  \leq C  \big(\|f\|_{L^2(\o_T)} + \|\nabla \psi\|_{L^p(\o)}^{\frac p2}\big)\,
\end{align}
for $k \in \mathbb N$.
Therefore, there exist
 an $\rn$--valued  function $U \in L^2(\o_T)$ and an $\mathbb R^{n\times n}$--valued function $V \in L^2(\o_T)$ such that
\begin{equation}\label{par25k}
|\nabla u_{k}|^{p-2} \nabla u_k\rightharpoonup U \quad \hbox{weakly in $L^2(\o_T)$}
\end{equation}
and
\begin{equation}\label{par26k}
\nabla (|\nabla u_{k}|^{p-2} \nabla u_{k})\rightharpoonup V \quad \hbox{weakly in $L^2(\o_T)$,}
\end{equation}
up to subsequences.
Thereby, $U \in L^2((0, T); W^{1,2}(\o))$ and 
\begin{equation}\label{par38}
V= \nabla U\,.
\end{equation}
Owing to \eqref{june0} and \eqref{par25k},
$U = |\nabla u|^{p-2} \nabla u$, whence
\begin{equation}\label{june6k}
|\nabla u_k|^{p-2} \nabla u_k  \rightharpoonup |\nabla u|^{p-2} \nabla u \quad \hbox{weakly in $L^2((0,T); W^{1,2}(\o))$.}
\end{equation}
 Hence, inequality \eqref{secondcap2} follows,   via \eqref{june8}  and \eqref{june6k}.
\qed

\medskip
\par\noindent
{\bf Proof of Theorem \ref{secondconvex}}. The proof is completely analogous to that of Theorem \ref{secondcap}, save that   properties \eqref{aug25} and \eqref{aug26} now follow from \cite[Theorem 2.3]{cmsecond}, that holds for any bounded open convex set.
 \qed

\medskip 
\par
We conclude with the proofs of Theorems \ref{secondcapsys} and \ref{secondconvexsys} for systems. They  require some minor variant with respect to those of the corresponding results for equations.

\medskip
\par\noindent
{\bf Proof of Theorem \ref{secondcapsys}}.
Since we are assuming that $p\geq 2$,  
$$f \in L^2(\o_T) \to L^{p'}(\o_T) \to  L^{p'}((0,T); W^{-1,p'}(\o)) \quad \hbox{and} \quad \psi \in W^{1,p}_0(\o) \to L^2(\o)\,.$$
 Therefore, there exists a unique weak solution $u \in C([0,T]; L^2(\o))\cap L^p((0,T); W^{1,p}_0(\o))$, with $u_t \in L^{p'}((0,T); W^{-1,p'}(\o))$,  to problem \eqref{eqdirichlet}. Moreover, owing to estimate  \eqref{solamaink} (that also holds when $N>1$) with $f_k$, $\psi_k$ and $u_k$ replaced by  $f$, $\psi$ and $u$, we have that
 $u_t \in L^2(\o _T)$ and there exists a constant $C$ such that
\begin{align}\label{solaN} 
\| u\|_{L^\infty((0,T); W^{1,p}_0(\o))}^{\frac p2}+ \|u_t\|_{L^2(\o_T)}& \leq C \big(\|f\|_{L^2(\o_T)} + \|\nabla \psi\|_{L^p(\o)}^{\frac p2}\big)  \,.
\end{align}
For  $t\in [0,T]$, define the functions $g^t: \o \to \mathbb R^N$ and $w^t : \o \to \mathbb R^N$ as
$$g^t(x) =f(x,t) - u_t(x,t) \quad  \hbox{and} \quad   w^t(x) = u(x,t) \quad \hbox{for $x \in \o$.}$$
%
The same argument as in the proof of  equation \eqref{june3bis} tells us that the function $w^t$ is the weak solution to the Dirichlet problem 
\begin{equation}\label{auxN}
\begin{cases}
- {\rm div} (|\nabla w^t|^{p-2}\nabla w^t ) = g^t  & {\rm in}\,\,\, \o \\
 w^t =0  &
{\rm on}\,\,\,
\partial \o \,
\end{cases}
\end{equation}
for a.e. $t\in (0, T)$. Observe that $w^t$ is actually the weak solution to problem \eqref{auxN} in the standard sense, since equation \eqref{june2}, and hence \eqref{june3bis} now hold for every $\varphi \in W^{1,p}_0(\o)$. Indeed, inasmuch as $p \geq 2$, one has that $L^p((0,T); W^{1,p}_0(\o)) \to L^2(\o)$, and since $\rho \varphi \in L^p((0,T); W^{1,p}_0(\o))$ if $\rho \in C^\infty _0(0, T)$ and $\varphi \in W^{1,p}_0(\o)$, one has that
$$  \langle u_t , \rho\, \varphi \rangle = \int _0^T \rho (t) \int _\o u_t (x,t) \varphi (x,t)\, dxdt\,.$$
Next,  fix any  $t\in (0, T)$ for which \eqref{auxN} holds, and let
$\{g^t_k\} \subset C^\infty_0(\o)$ be a sequence such that $g^t_k \to g^t$ in $L^2(\o)$ and $\|g^t_k\|_{L^2(\o)} \leq 2 \|g^t\|_{L^2(\o)}$ for $k \in \mathbb N$. For each $k$, let $w^t_k$ be the weak solution to the problem
\begin{equation}\label{auxNk}
\begin{cases}
- {\rm div} (|\nabla w^t_k|^{p-2}\nabla w^t_k ) = g^t_k  & {\rm in}\,\,\, \o \\
 w^t_k =0  &
{\rm on}\,\,\,
\partial \o \,.
\end{cases}
\end{equation}
The use of $w^t_k$ as a test function in the weak formulation of problem \eqref{auxNk} and of the H\"older and  the Sobolev inequalities tells  us that 
\begin{equation}\label{auxN2}
\|w^t_k\|_{L^2(\o)} \leq c \|g^t\|_{L^2(\o)}
\end{equation}
for $k \in\mathbb N$, for some constant $c=c(n,N,p, |\o|)$. On the other hand, 
choosing  $w^t_k - w^t$ as a test function in the weak formulations of problems \eqref{auxN} and \eqref{auxNk}, subtracting the resultant equations, and observing that 
$$c |\xi - \eta|^p \leq (|\xi|^{p-2} \xi - |\eta |^{p-2}\eta) \cdot (\xi - \eta) \quad \hbox{for $\xi, \eta \in \mathbb R^{N\times n}$,}$$
for some positive constant $c=c(n,N,p)$ yield
\begin{align}\label{auxN1}
c \int _\o |\nabla w^t_k - \nabla w^t|^p\, dx&  \leq 
\int_\o(|\nabla w^t_k|^{p-2}\nabla w^t_k - |\nabla w^t|^{p-2}\nabla w^t)\cdot (\nabla w^t_k - \nabla w^t)\, dx \\ \nonumber & = \int _\o (g^t_k -g^t) \cdot (w^t_k-w^t) \, dx  \leq 
\|g^t_k -g^t\|_{L^2(\o)} \|w^t_k-w^t\|_{L^2(\o)} \\ \nonumber &\leq  \|g^t_k -g^t\|_{L^2(\o)} (\|w^t_k\|_{L^2(\o)}+ \|w^t\|_{L^2(\o)})
\end{align}
for $k \in\mathbb N$. Hence, by inequality \eqref{auxN2}, $w^t_k \to w^t$ in $W^{1,p}_0(\o)$ as $k \to \infty$. As a conequence, for a.e. $t\in (0,T)$ the function $w^t$ is an approximable solution to problem \eqref{auxN}. An application of
  \cite[Theorem 2.1]{cmsecondsys} then entails that
$ |\nabla u(\cdot , t)|^{p-2} \nabla u(\cdot , t)= |\nabla w^t(\cdot)|^{p-2}\nabla w^t(\cdot) \in W^{1,2}(\o)$ and 
\begin{equation}\label{june5N}
\||\nabla u(\cdot , t)|^{p-2} \nabla u(\cdot , t)\|_{W^{1,2}(\o)}^2 \leq C \int_\o f(x,t)^2 + u_t(x,t)^2\, dx \quad \hbox{for a.e. $t \in (0,T)$,}
\end{equation}
some constant  $C=C(n,N,p, \o)$. Equations \eqref{secondcapsys1} and \eqref{secondcapsys2} follow on 
integrating inequality \eqref{june5N} with respect to $t$ over $(0,T)$, and exploiting inequality \eqref{solaN}. \qed

\medskip
\par\noindent
{\bf Proof of Theorem \ref{secondconvexsys}}. The proof is the same as that of Theorem \ref{secondcapsys},  inequality \eqref{june5N}  being  a consequence of \cite[Theorem 2.6]{cmsecondsys} under the assumption that $\o$ is a bounded convex set.
 \qed

\color{black}

\bigskip
\par\noindent
{\bf Acknowledgment}. This research was initiated during a visit of the first-named author at the Institut Mittag-Leffler in August 2017. He wishes to thank the Director and the Staff of the Institut for their support and hospitality.

\color{black}
\section*{Compliance with Ethical Standards}

\smallskip
\par\noindent 
{\bf Funding}. This research was partly funded by:   
\\ (i) Research Project of the
Italian Ministry of University and Research (MIUR) Prin 2012 \lq\lq Elliptic and
parabolic partial differential equations: geometric aspects, related
inequalities, and applications"  (grant number 2012TC7588);    
\\ (ii) GNAMPA   of the Italian INdAM - National Institute of High Mathematics (grant number not available);   

\smallskip
\par\noindent
{\bf Conflict of Interest}. The authors declare that they have no conflict of interest.

\end{document}